\documentclass{ifacconf}

\usepackage{graphicx}      % include this line if your document contains figures
\usepackage{natbib}        % required for bibliography
\usepackage{amsmath}
\usepackage{amssymb} 

\newtheorem{defin}{Definition}

\begin{document}
\begin{frontmatter}

\title{Novel Variants of Diffusive Representation of Fractional Integrals: Construction and Numerical Computation\thanksref{footnoteinfo}} 
\thanks[footnoteinfo]{This work has been supported by 
	the German Federal Ministry for Education and Research (BMBF)
	under Grant No.\ 05M22WHA.}

\author[First]{Renu Chaudhary} 
\author[First]{Kai Diethelm}

\address[First]{FANG, Technical University of Applied Sciences Würzburg-Schweinfurt, Ignaz-Schön-Str. 11, 97421 Schweinfurt, Germany (Email: \{renu.chaudhary, kai.diethelm\}@thws.de)}

\begin{abstract}                % Abstract of not more than 250 words.
In this paper, we revisit the diffusive representations of fractional integrals established in \cite{diethelm2023diffusive} to explore novel variants of such representations which provide highly efficient numerical algorithms for the approximate numerical evaluation of fractional integrals.
\end{abstract}


\end{frontmatter}
%===============================================================================
\section{Introduction}
Fractional calculus, nowadays an established mathematical discipline, plays a pivotal role in understanding complex physical and mathematical phenomena. The application of fractional calculus extends to various fields such as physics, engineering, finance, and signal processing, among others. For a comprehensive introduction and exploration of its applications, we refer to \cite{Podlubny}, \cite{Kilbas} and \cite{Diethelmbook}. Fractional integrals, a fundamental concept within this domain, have been extensively studied and applied for solving real-world problems.

In certain scenarios, discovering analytical solutions for problems involving fractional order operators can be challenging or even practically impossible. Hence a large number of different approaches have been proposed for numerical approximation of fractional order differential and integral operators. Due to the intrinsic non-local nature of fractional order differential and integral operators, achieving their numerical evaluation in the conventional representation presents a significantly heightened computational challenge compared to assessing their integer order counterparts. This challenge pertains to both run-time efficiency and memory demands, particularly when evaluating at multiple points. For instance, conventional algorithms such as fractional linear multi-step methods in \cite{multistep} and Adams method in \cite{predictor, erroranalysis} may pose computational challenges due to their reliance on numerical approximations incurring a computational cost of  $\mathcal{O}(P^2)$ and a storage cost of $\mathcal{O}(P)$, where $P$ represents the number of evaluation points.

In this case, employing diffusive representations of fractional operators presents a dynamic and sustainable solution to address the issue at hand. Utilizing diffusive representations of fractional differential and integral operators offers a practical avenue for crafting highly efficient numerical algorithms aimed at their approximate assessment with less memory footprint $\mathcal{O}(1)$ and arithmetic complexity $\mathcal{O}(P)$, which aligns with the techniques commonly employed in solving problems with integer order operators. For more comprehensive information on diffusive representation, please refer to \cite{montseny1998diffusive}, \cite{diethelm2022new} and \cite{diethelm2021new}. 

In this research, we delve into the realm of fractional integrals, focusing on their diffusive representations. These representations serve as a bridge between theoretical understanding and practical applications. Building upon the foundational work by \cite{diethelm2023diffusive}, we present novel variants of these diffusive representations, pushing the boundaries of efficiency and accuracy in numerical computation of fractional integrals.

Our objective is to enhance the existing diffusive representations, offering improved algorithms for approximating fractional integrals numerically. By doing so, we aim to provide researchers and practitioners with powerful tools to handle the intricate calculations associated with fractional integrals. These enhanced numerical methods open doors to solving complex problems with a higher degree of precision and efficiency, thus advancing the applicability and impact of fractional calculus in diverse domains. In the subsequent sections, we elucidate the existing diffusive representations, introduce our novel variant, and showcase how these variants pave the way for highly efficient numerical algorithms. Additionally, we outline potential future directions in this exciting area of research.

\section{Mathematical Background}
 
Our aim is to explore novel variants of diffusive representation for Riemann-Liouville fractional integrals of order $\alpha\in (0,1)$ for functions $f\in C[a,b]$ (where $a, b \in \mathbb R$ such that $a<b$) with a starting point at $a$ represented as
\begin{equation}\label{RLIntegral}
I^{\alpha}_a f(t) = \frac{1}{\Gamma(\alpha)}\int_a^t(t-\tau)^{\alpha-1}f(\tau) \, \mathrm d\tau.
\end{equation}

As discussed in \cite{diethelm2023diffusive}, the diffusive representation of (\ref{RLIntegral}) can be formulated as an integral involving an auxiliary bivariate function such that the integration is performed with respect to a single variable, within a fixed range. This characteristic allows for numerical computation using a fixed quadrature formula whose arithmetic complexity does not increase over time. Additionally, the auxiliary function can be determined as the unique solution to an initial value problem for a first-order differential equation, obviating the need to consider memory effects.
 
\begin{defin}
Let $\Omega$ be a non-empty open interval. A strictly monotonically increasing function $\psi: \Omega \rightarrow (0,\infty)$ is called an admissible transformation if it possesses the following two properties:
\begin{enumerate}
  \item[(a)] $\psi\in C^1(\Omega)$,
  \item[(b)] $\lim_{r\rightarrow \inf \Omega}\psi(r)=0$ and $\lim_{r\rightarrow \sup \Omega}\psi(r)=+\infty$. 
\end{enumerate}
\end{defin}

\begin{thm} For a given admissible transformation  $\psi: \Omega \rightarrow (0,\infty)$, the Riemann-Liouville integral of order $\alpha\in (0,1)$ for the function $f\in C[a,b]$, where $a, b \in \mathbb R$, can be articulated in the following diffusive representation
\begin{equation}\label{RLDiffusiveIntegral}
I^{\alpha}_af(t)=\int_{\Omega}\phi(t,r) \, \mathrm dr,
\end{equation}
where
\begin{equation}\label{phi}
\phi(t,r)=c_{\alpha}\psi^{\prime}(r)(\psi(r))^{-\alpha}\int_a^t \mathrm e^{-(t-\tau) \psi(r)} f(\tau) \, \mathrm d\tau,
\end{equation}
with constant $c_{\alpha}=\pi^{-1} \sin\pi \alpha$.

Moreover, for a fixed value of  $r \in \Omega$, the function $\phi(\cdot, r)$ is characterized as the unique solution to the following initial value problem for first-order differential equation on the interval $[a, b]$
\begin{eqnarray}\label{DE}
\begin{split}
\frac{\partial \phi(t,r)}{\partial t}&=-\psi(r) \phi(t,r) + c_{\alpha}\psi^{\prime}(r)(\psi(r))^{-\alpha}f(t),\\
\phi(a,r)&=0.
\end{split}
\end{eqnarray}
\end{thm}
\begin{pf}
The proof of the theorem can be straightforwardly deduced from the proof presented in \citet[Theorems 1 and 2]{diethelm2023diffusive}, specifically for the case $n=1$. 
\qed
\end{pf}

In literature, numerous approaches have been devised to articulate the fractional integral of a provided function. For instance, in the work by \cite{RebeccaLi}, the author computed the fractional integral by employing an integral representation of the convolution kernel. She proceeded to devise an effective quadrature for this integral representation and incorporated it into a fast time-stepping technique. In \cite{beylkin}, the authors applied the Poisson's summation formula to discretize the integral representations of power functions to obtain their approximations by exponential functions. Later on, in \cite{mclean2018exponential}, the author extended the idea of \cite{beylkin} by introducing an alternative integral representation of the power function. In the study \cite{baffet}, the Gauss-Jacobi quadrature rule was utilized to approximate the kernel of the fractional integral through a linear combination of exponentials.

The main objective of this paper is to devise a diffusive representation that facilitates more efficient numerical handling of the integral (\ref{RLDiffusiveIntegral}) than traditional approaches. For this, we have examined an important special case of admissible transformation $\psi(\cdot)$ such that the corresponding bivariate function $\phi(t,\cdot)$ demonstrates an exponential decay at the endpoints of the domain $\Omega$. This behavior holds significant value from an approximation theoretic perspective which allows us to use truncated trapezoidal scheme or the Gauss-Laguerre numerical quadrature scheme for integrating the function $\phi$. Moreover, it is imperative to have knowledge of the integrand $\phi(t,\cdot)$, which we have achieved by numerically solving the associated first-order differential equation (\ref{DE}) through Euler's Backward method or the trapezoidal method.

\section{Novel Variants of Diffusive Representation of Fractional Integrals}
In this section, we discuss a diffusive representation to the integral (\ref{RLDiffusiveIntegral}) using the well-known exponential function as an admissible transformation. We prove that in this new representation the integrand $\phi(t,r)$ showcases exponential decay as $r$ tends towards $\pm\infty$. This characteristic, coupled with a smoothness result, allows a highly efficient numerical integration.

\begin{thm}\label{mainthm1}
	Consider a function $f\in C[a,b]$, where $a, b \in \mathbb R$ such that $a<b$, and $\alpha\in (0,1)$. 
	Let us define the function $\phi(t,r)$ as
	\begin{equation}\label{ExDiffusiveIntegral2}
		\phi(t,r)=c_{\alpha} \mathrm e^{(1-\alpha)r} \int_a^t \mathrm e^{-(t-\tau) \mathrm e^r} f(\tau) \, \mathrm d\tau.
	\end{equation}
	for all $r\in \mathbb{R}$ and $t\in [a,b]$. The following properties can then be observed:
	\begin{enumerate}
	\item[(i)] The function $\phi(\cdot, r)$ is characterized as the unique solution to the 
		following initial value problem for first-order differential equation on the interval $[a, b]$
		\begin{eqnarray}\label{DE10}
			\begin{split}
				\frac{\partial \phi(t,r)}{\partial t}&=- \mathrm e^r \phi(t,r) + c_{\alpha} \mathrm e^{(1-\alpha)r}f(t),\\
				\phi(a,r)&=0.
			\end{split}
		\end{eqnarray}
	\item[(ii)] For any $t\in [a,b]$, 
		\begin{equation}\label{ExDiffusiveIntegral1}
			I^{\alpha}_af(t)=\int_{-\infty}^{\infty} \phi(t,r) \, \mathrm dr.
		\end{equation}
	\item[(iii)] For any $t\in [a,b]$, the integrand $\phi(t,\cdot)\in C^{\infty}(\mathbb{R})$.
 	\item[(iv)] For any $t$ belonging to the interval $[a, b]$, there exists a real constant $C$ such that:
		\begin{equation}\label{asymp 1}
			 |\phi(t,r)|\leq C \mathrm e^{-\alpha r}\quad \textit{as}\,\,\, r \rightarrow \infty,
		\end{equation}
		and 
		\begin{equation}\label{asymp 2}
			|\phi(t,r)|\leq C \mathrm e^{(1-\alpha) r} \quad \textit{as}\,\,r\rightarrow -\infty.
		\end{equation}
	\end{enumerate}
\end{thm}

\begin{pf}
By selecting $\psi(r)= \mathrm e^r$ with $\Omega=(-\infty,\infty)$ in (\ref{RLDiffusiveIntegral}) and (\ref{phi}), the proof of the Theorem can be directly deduced from the results presented in \cite{diethelm2023diffusive}.
\qed
\end{pf}

As emphasized in \cite{diethelm2023diffusive}, the accuracy of the numerical approximation for the integral (\ref{ExDiffusiveIntegral1}) is significantly influenced by two key aspects of the integrand. The first crucial aspect is the smoothness of the integrand, while the second pertains to the asymptotic behavior of the integrand as the integration variable approaches the integration limit. In the context of a diffusive representation, these two critical aspects of the integrand are profoundly influenced by the choice of the admissible transformation, which in this case is given by the exponential function.

While the general representation \eqref{RLDiffusiveIntegral} and its special case \eqref{ExDiffusiveIntegral2} 
have proven to be very useful in analytical considerations \citep{montseny1998diffusive}, 
it has been observed \citep{mclean2018exponential} that there are certain advantages from the numerical approximation
perspective in handling the operators in a slightly different way. Specifically, for $t \ge a+h$ with some $h > 0$ it is useful to split up the hereditary integral from
\eqref{RLIntegral} into a local part $L_{a, h}^\alpha f(t)$ and a history part $H_{a, h}^\alpha f(t)$ 
according to
\[
	I_a^\alpha f(t) = L_{a, h}^\alpha f(t) + H_{a, h}^\alpha f(t)
\]
where
\begin{equation}
	\label{eq:def-l}
	 L_{a, h}^\alpha f(t) = \frac{1}{\Gamma(\alpha)} \int_{t-h}^t (t-\tau)^{\alpha-1} f(\tau) \, \mathrm d\tau
\end{equation}
and
\begin{equation}
	\label{eq:def-l}
	 H_{a, h}^\alpha f(t) = \frac{1}{\Gamma(\alpha)} \int_a^{t-h} (t-\tau)^{\alpha-1} f(\tau) \, \mathrm d\tau.
\end{equation}
The local part $L_{a,h}^\alpha f$ can usually be handled directly, but the history part $H_{a,h}^\alpha f(t)$ benefits
from a representation that is analog to the diffusive form for the full Riemann-Liouville integral discussed above.
Specifically, the following relationships hold.

\begin{thm}
	\label{thm:diff-rep-hist}
	Let $f \in C[a, b]$ with some real numbers $a < b$, and
	let $\alpha \in (0,1)$.
	Moreover, assume that $\psi : \Omega \to (0, \infty)$ 
	is an admissible transformation and that $0 < h < b-a$.
	Then, for every $t \in [a + h, b]$, the history part of the Riemann-Liouville integral
	of order $\alpha$ of the function $f$ can be expressed in the
	form of the diffusive representation
	\begin{subequations}
		\label{eq:diff-rep-hist-complete}
	\begin{equation}
		\label{eq:diff-rep-hist}
		H_{a,h}^\alpha f(t) = \int_\Omega \mu(t, h, r) \, \mathrm d r
	\end{equation}
	with
	\begin{align}
		\label{eq:mu1}
		\mu(t, h, r) 
		&= c_\alpha \psi'(r) (\psi(r))^{-\alpha} \\
		& \nonumber 
			\qquad \times
			\int_a^{t-h} \exp(- (t - \tau) 
						\psi(r)) f(\tau) \, \mathrm d \tau.
	\end{align}
	\end{subequations}
	The function $\mu$ from \eqref{eq:mu1} has the following properties:
	\begin{enumerate}
	\item[(a)] For any $r \in \Omega$, 
		the function $\mu(\cdot, h, r)$ is the unique solution 
		on the interval $[a+h,b]$ to the first order differential equation
		\begin{subequations}
			\label{eq:ivp-1}
		\begin{eqnarray}
			\label{eq:ode-1}
			\frac{\partial \mu}{\partial t} (t, h, r)  & = & - \psi(r) \mu(t, h, r) \\
			\nonumber
			&& {} + c_\alpha \psi'(r) (\psi(r))^{-\alpha} \mathrm e^{-h \psi(r)} f(t - h)
		\end{eqnarray}
		subject to the initial condition
		\begin{equation}
			\label{eq:ic-1}
			\mu(a+h, h, r) = 0.
		\end{equation}
		\end{subequations}
	\item[(b)] If the transformation function $\psi$ satisfies $\psi \in C^k(\Omega)$
		with some $k \in \mathbb N$ then, for any $t \in [a+h, b]$,
		we have $\mu(t, h, \cdot) \in C^{k-1}(\Omega)$.
	\item[(c)] If the transformation function $\psi$ satisfies 
		$\psi \in C^\infty(\Omega)$
		then, for any $t \in [a+h, b]$,
		we have $\mu(t, h, \cdot) \in C^\infty(\Omega)$.
	\item[(d)] If $f \in C^\ell[a, b]$ with some $\ell \in \mathbb N_0$
		then, for any $r \in \Omega$, $\mu(\cdot, h, r) \in C^{\ell+1}[a+h, b]$.
	\item[(e)] For any fixed $t \in [a+h, b]$, there exist some 
		constants $C_1, C_2 > 0$ such that
		\begin{align}
			\label{eq:mu-asymp-sup}
			|\mu(t, h, r)| &\le C_1 \psi'(r) (\psi(r))^{-\alpha-1} \mathrm e^{-h \psi(r)} \\
			& \nonumber
			\qquad \mbox{ for } r \to \sup \Omega \\
		\intertext{and}
			\label{eq:mu-asymp-inf}
			|\mu(t, h, r)| & \le C_2 \psi'(r) (\psi(r))^{-\alpha}
			\quad \mbox{ for } r \to \inf \Omega.
		\end{align}
	\end{enumerate}
\end{thm}

\begin{rem}
	With respect to Theorem \ref{thm:diff-rep-hist} (e), we note that the asymptotic estimate
	\eqref{eq:mu-asymp-inf} for $r \to \inf \Omega$ shows the same rate of decay for this
	history term as the corresponding result for the full Riemann-Liouville integral
	\cite[part 2 of Theorem 5]{diethelm2023diffusive}. Our estimate \eqref{eq:mu-asymp-sup} 
	for $r \to \sup \Omega$, however, shows a much faster decay than the corresponding
	result for the full integral from \citet[part 1 of Theorem 5]{diethelm2023diffusive}.
	To the best of our knowledge, this fact has not been fully exploited in the construction
	of numerical methods yet.
\end{rem}

\begin{pf}
	The proof of the representation \eqref{eq:diff-rep-hist-complete} for the history part of the
	Riemann-Liouville integral is completely analog to the proof of the corresponding relation
	for the complete Riemann-Liouville integral provided in \citet[Section V.A]{diethelm2023diffusive}.
	
	To show the properties of the function $\mu$, we proceed as follows:
	The differential equation in (a) can be verified by a direct differentiation of the function $\mu$
	with respect to the first argument; the initial condition is an immediate consequence 
	of the representation \eqref{eq:mu1}.

	Items (b) and (c) also follow from \eqref{eq:mu1}.
	
	For part (d), we use the notation $M(t) = \mu(t, h, r)$ for some arbitrary but fixed $r \in \Omega$
	and $h > 0$. Then, the differential equation \eqref{eq:ode-1} can be rewritten as
	$M'(t) = G(t, M(t))$ where 
	$G(t, z) = -\psi(r) z + c_\alpha \psi'(r) (\psi(r))^{-\alpha} \mathrm e^{-h \psi(r)} f(t - h)$.
	Clearly, under the assumption of (d), $G \in C^\ell([a+h, b] \times \mathbb R)$.
	and therefore a standard argument from the theory of differential equations 
	\cite[Chapter~1, Theorem 1.2]{CL} yields the claim.
	
	Finally, for (e) we see that \eqref{eq:mu1} implies
	\[
		| \mu(t, h, r) | \leq c_\alpha \sup_{\tau \in [a,b]} |f(\tau)| \cdot \psi'(r) (\psi(r))^{-\alpha} J(r)
	\]
	with 
	\[
		J(r) = \int_a^{t-h} \mathrm e^{- (t - \tau) \psi(r)} \, \mathrm d \tau
			 = \frac 1 {\psi(r)} \int_{h \psi(r)}^{(t-a) \psi(r)} \mathrm e^{-\rho} \, \mathrm d \rho.
	\]
	Thus, when $r \to \inf \Omega$ and hence $\psi(r) \to 0$, then (since the integrand is always
	less than $1$) 
	\[
		J(r) \le \frac 1 {\psi(r)} \int_{h \psi(r)}^{(t-a) \psi(r)} \, \mathrm d \rho  = t - a - h,
	\]
	and so \eqref{eq:mu-asymp-inf} follows. For $r \to \sup \Omega$ on the
	other hand, we find that
	\[
		J(r) \le \frac 1 {\psi(r)} \int_{h \psi(r)}^{\infty} \mathrm e^{-\rho} \, \mathrm d \rho 
		= \frac 1 {\psi(r)} \mathrm e^{-h \psi(r)}
	\]
	which proves \eqref{eq:mu-asymp-sup}.
	\qed
\end{pf}

The result of Theorem \ref{thm:diff-rep-hist} can be extended to the case
of general $\alpha > 0$, $\alpha \notin \mathbb N$. We shall discuss this extension
in a separate paper in the future.

\subsection{Kernel Function Approximation by Exponential Sums}\label{expsumapp1}

Now we derive the exponential sum approximation to the kernel function of the integral represented by (\ref{ExDiffusiveIntegral1}). Here we utilize the well established trapezoidal rule to discretize the integral. This rule provides an explicit discretization of the integral over the unbounded domain, expressed as a sum of exponentials. Our integrand in (\ref{ExDiffusiveIntegral2}) exhibits the crucial property of exponential decay, as outlined in Theorem \ref{mainthm1}(iv). This property permits the effective utilization of the trapezoidal rule. For a deeper understanding of the exponential sum approximation, one may refer to the work of \cite{beylkin} and \cite{mclean2018exponential}.

Utilizing the proof presented in \citet[Theorem 1]{diethelm2023diffusive}, we can express the value of $I^{\alpha}_af(t)$ as
\begin{equation*}
I^{\alpha}_af(t) = c_{\alpha} \int_a^t \int_{0}^{\infty} \bigg(\frac{u}{t-\tau}\bigg)^{(1-\alpha)} \mathrm e^{-u}u^{-1} f(\tau) \, \mathrm d u \, \mathrm d\tau.
\end{equation*}
By making the substitution $u=(t-\tau) \mathrm e^r$, we arrive at
\begin{equation*}
I^{\alpha}_af(t)=c_{\alpha}\int_a^t\int_{-\infty}^{\infty} \mathrm e^{(1-\alpha)r}\exp[-(t-\tau) \mathrm e^r]f(\tau)\, \mathrm d r\, \mathrm d\tau
\end{equation*}
which can be further written as
\begin{eqnarray*}
I^{\alpha}_af(t)= c_{\alpha}\int_a^tK(t-\tau)f(\tau)\, \mathrm d\tau,
\end{eqnarray*}
where the kernel $K(t)$ is given by
 \begin{equation}
K(t)=\int_{-\infty}^{\infty} \mathrm e^{(1-\alpha)r} \mathrm e^{-t \mathrm e^r}\, \mathrm d r\quad \text{for}\quad t\in [0,b-a].
\end{equation}
Now applying the trapezoidal rule to discretize the infinite integral with step size $h>0$, we obtain
\begin{eqnarray*}
K(t)&\approx&h\sum_{n=-\infty}^{\infty} \mathrm e^{(1-\alpha)nh} \mathrm e^{-t \mathrm e^{nh}}\\
&=&\sum_{n=-\infty}^{\infty} \tilde w_n \mathrm e^{-\tilde \beta_nt}
\end{eqnarray*}
where 
\begin{equation}\label{expoandweights}
  \tilde w_n = h \mathrm e^{(1-\alpha)n h}\quad \text{ and } \quad \tilde \beta_n = \mathrm e^{n h}.
\end{equation}
If now $t$ is restricted to a compact interval $[\delta,b-a]$ with $0<\delta<b-a<\infty$, we can proceed in this approach to obtain a finite exponential sum approximation
\begin{eqnarray}\label{finitesumappr}
	K(t)\approx\sum_{n=-M}^{N} \tilde w_n \mathrm e^{-\tilde \beta_n t} \quad \text{ for }\,\,t\in [\delta,b-a]
\end{eqnarray}
for the kernel with suitably chosen $M, N \in \mathbb N$.

\citet{beylkin} discussed an innovative reduction algorithm based on Prony's method, specifically designed for cases with excessively many terms and small exponents in the exponential sum approximation to a kernel function. This suggests that employing Prony's method allows for a substantial reduction in the number of terms in the exponential sum in (\ref{finitesumappr}) without compromising accuracy. Consequently, numerical algorithms designed using this exponential sum approximation for approximating the fractional integral demonstrate markedly enhanced efficiency when compared to existing numerical methods.

\subsection{Truncation Error Estimation of the Kernel Function}
When conducting practical computations, particularly in numerical analysis and real-world applications, finite computational resources necessitate the use of truncated sums. In such cases, we need to approximate the truncation error for the exponential sum defined in (\ref{finitesumappr}). In this section, we assess the estimation of the two components of the truncation error:

\begin{thm}
If $t \in [0, b-a]$, the exponents and weights are given by (\ref{expoandweights}) and
\begin{align*}
t \mathrm e^{Nh} \geq 1-\alpha \geq t \mathrm e^{-Mh},
\end{align*}
then 
\begin{equation*}
\sum_{n=N+1}^{\infty} \tilde w_n \mathrm e^{- \tilde \beta_nt}\leq \frac{1}{t^{(1-\alpha)}} \Gamma(1-\alpha,t \mathrm e^{Nh})
\end{equation*}
and 
\begin{equation*}
\sum_{n=-\infty}^{-M-1} \tilde w_n \mathrm e^{- \tilde \beta_n t}
\leq \frac{1}{t^{(1-\alpha)}} \bigg(\Gamma(1-\alpha)-\Gamma(1-\alpha,t \mathrm e^{-Mh})\bigg)
\end{equation*}
where $\Gamma(\cdot, \cdot)$ denotes the upper incomplete Gamma function.
\end{thm}
\begin{pf}
The proof directly follows from \citet[Theorem 2]{mclean2018exponential} with appropriate modifications to the conditions and parameters. \qed
\end{pf}

\section{Efficient Numerical Algorithms for Approximate Evaluation of Fractional Integral}\label{Numerical Algorithms}
In this section, we establish two distinct efficient numerical algorithms for the approximate evaluation of the fractional integral represented in (\ref{ExDiffusiveIntegral1}). The first approach involves utilizing the exponential sum approximation to the kernel function, as defined in (\ref{finitesumappr}). The second approach employs the Gauss-Laguerre quadrature formula.

\subsection{Numerical Algorithm for Approximating Fractional Integral using Exponential Sum Approximation of the Kernel Function}
For the numerical computation of $I_a^{\alpha}f(t)$ at the points of the grid $a=t_0<t_1<\cdots <t_{P}=b$, let $f(t_n)\approx F^n$ and define a piecewise constant interpolant $\tilde{F}(t)=F^n$ for $t\in J_n=(t_{n-1},t_n)$, then
\begin{equation}
  I^{\alpha}_af(t_n) =c_{\alpha}\int_a^{t_n}K(t_n-\tau)f(\tau)\, \mathrm d \tau\approx\sum_{j=1}^{n}z_{nj}F^j
\end{equation}
where 
\begin{equation}
z_{nj}=c_{\alpha}\int_{J_j}K(t_n-\tau)\, \mathrm d\tau.
\end{equation}

The computational complexity for calculating this sum within the range of $1\leq n\leq P$ scales proportionally to $P^2$. Such quadratic growth can pose challenges, particularly in scenarios where each $F^j$ represents a large vector rather than a scalar. In particular, storing $F^j$ in active memory for all time levels $j$ may not be feasible. 

Here we define an efficient algorithm to avoid these challenges based on the exponential sum approximation to the kernel $K$ such as 
\begin{equation}
K(t)\approx \sum_{l=1}^{\Lambda}w_l \mathrm e^{\beta_l t} \quad \text{for}\quad \delta\leq t\leq b-a, 
\end{equation} 
provided that the moderate number of terms $\Lambda$ can achieve adequate accuracy for a choice of $\delta>0$ that is smaller than the time step $\Delta t_n=t_n-t_{n-1}$ for all $n$. Certainly, if $\Delta t_n\geq \delta$ then $\delta\leq t_n-\tau \leq b-a$ for $0\leq \tau \leq t_{n-1}$; therefore
\begin{eqnarray*}
\sum_{j=1}^{n-1}z_{nj}F^j &=&c_{\alpha}\int_{0}^{t_{n-1}}K(t_n-\tau)\tilde{F}(\tau)\, \mathrm d\tau\\
&\approx& c_{\alpha}\int_{0}^{t_{n-1}}\sum_{l=1}^{\Lambda}w_l \mathrm e^{\beta_l (t_n-\tau)}\tilde{F}(\tau)\, \mathrm d\tau\\
&=& \sum_{l=1}^{\Lambda} \Phi_l^n
\end{eqnarray*}
where 
\begin{equation*}
\Phi_l^n =c_{\alpha}w_l \int_{0}^{t_{n-1}} \mathrm e^{\beta_l (t_n-\tau)}\tilde{F}(\tau)\, \mathrm d\tau
= \sum_{j=1}^{n-1}K_{lnj}F^j
\end{equation*}
with 
\begin{equation}
K_{lnj}= c_{\alpha}w_l\int_{J_j} \mathrm e^{\beta_l (t_n-\tau)}\, \mathrm d\tau.
\end{equation}
Hence we express the approximated value of the integral
\begin{equation}\label{Trapappro1}
I^{\alpha}_a\tilde{F}(t_n)\approx z_{nn}F^n+ \sum_{l=1}^{\Lambda} \Phi_l^n,
\end{equation}
and utilizing the recursive relation
\begin{align}
\Phi_l^1 &= 0, \\
\label{Trapappro2}
\Phi_l^n &= K_{ln,n-1} F^{n-1} + \mathrm e^{\beta l \Delta t_n}\Phi_l^{n-1} \quad (n = 2, 3, \ldots, P),
\end{align}
we obtain an evaluation scheme that, at time $t_n$, only needs to exploit information about the immediately preceding time $t_{n-1}$
but not about any earlier point in history. 

Thus we can efficiently compute $I^{\alpha}_a\tilde{F}(t_n)$ for $n = 1, 2, \ldots, P$ with a satisfactory level of accuracy using a number of operations proportional to $\Lambda \cdot P$. This approach yields substantial computational savings compared to the common $\mathcal O(P^2)$ operation count for traditional approaches, particularly when $\Lambda \ll P$. Moreover, we have the flexibility to overwrite $\Phi_l^{n-1}$ with $\Phi_l^{n}$, and $F^{n-1}$ with $F^{n}$, effectively reducing the active storage requirement from being proportional to $P$ to being proportional to $\Lambda$.

\subsection{Numerical Algorithm for Approximating Fractional Integral using Gauss-Laguerre Quadrature Formula}\label{GLQuadrature}

Next, we approximate the integral (\ref{ExDiffusiveIntegral1}) using the $\Lambda$-point Gauss-Laguerre quadrature formula. This quadrature method is an open-type Gaussian quadrature rule designed for evaluating integrals, specifically with the weight function $\mathrm e^{-u}$, over the interval $[0,\infty)$. We write
\begin{eqnarray*}
I^{\alpha}_af(t) &=& \int_{-\infty}^{\infty}\phi(t,r)\, \mathrm dr\\
&= &\int_{-\infty}^{0}\phi(t,r)dr +\int_{0}^{\infty}\phi(t,r)\, \mathrm dr\\
&=&J_1+J_2.
\end{eqnarray*}
Substituting $r=\frac{-s}{1-\alpha}$ in $J_1$ and $r=\frac{s}{\alpha}$ in $J_2$, we obtain
\begin{equation*}
\begin{split}
 I^{\alpha}_a f(t) = {}& \frac{1}{1-\alpha}\int_{0}^{\infty} \mathrm e^{-s} \mathrm e^{s}\phi(t,-s/(1-\alpha))\, \mathrm ds\\
&+\frac{1}{\alpha}\int_{0}^{\infty} \mathrm e^{-s} \mathrm e^{s}\phi(t,s/\alpha)\, \mathrm ds.
\end{split}
\end{equation*}
Thus, employing
\begin{eqnarray*}
\hat{\phi}(t,s) =& \mathrm e^{s} \bigg(\frac{1}{1-\alpha}\phi(t,-s/(1-\alpha))+\frac{1}{\alpha}\phi(t,s/\alpha)\bigg),
\end{eqnarray*}
we obtain that 
\begin{eqnarray*}
I^{\alpha}_af(t) =\int_{0}^{\infty} \mathrm e^{-s}\hat{\phi}(t,s)\, \mathrm ds\approx  Q_\Lambda^{GLa}[\hat{\phi}(t,\cdot)],
\end{eqnarray*}
where the expression
\begin{eqnarray*}
Q_\Lambda^{GLa}[\hat{\phi}]=\sum_{l=1}^{\Lambda} w_l^{GLa} \hat{\phi}(x_l^{GLa})
\end{eqnarray*}
represents the $\Lambda$-point Gauss-Laguerre quadrature formula with weights $w_l^{GLa}$ and nodes $x_l^{GLa}$. 

To determine the Gauss-Laguerre nodes $x_l^{GLa}$ and the corresponding weights $w_l^{GLa}$, we need to find the zeros of the Laguerre polynomial $L_\Lambda$ of order $\Lambda$ and calculate the weights using a specific formula:
\begin{itemize}
   \item Solve $L_\Lambda(x_l^{GLa})=0$, for $l=1,2,\ldots, \Lambda$ to obtain the Gauss-Laguerre nodes  $x_l^{GLa}$.
   \item Use the provided formula to calculate the weights:
   \begin{eqnarray*}
	w_l^{GLa}=\frac{x_l^{GLa}}{[L_{\Lambda+1}(x_l)^2]}.
   \end{eqnarray*}
\end{itemize}
%The nodes $x_n^{GLa}$ with $n=1,2,\ldots N$, correspond to the zeros of the Laguerre polynomial $L_N$ of order $N$ and their associated weights are calculated as follows:
%\begin{eqnarray*}
%w_n^{GLa}=\frac{x_n^{GLa}}{[L_{N+1}(x_n)^2]}.
%\end{eqnarray*}
For a more comprehensive understanding of Gauss-Laguerre quadrature formulas, refer to \cite{davis2007methods}.

\begin{thm}
Under the assumptions of Theorem \ref{mainthm1}, 
\begin{equation*}
\lim_{\Lambda \rightarrow \infty} Q_\Lambda^{GLa}[\hat{\phi}(t,\cdot)] = I^{\alpha}_a f(t) 
\end{equation*}
for all $t\in [a,b]$.
\end{thm}
\begin{pf} By Theorem \ref{mainthm1}(iii), the function $\phi(t,\cdot)$ possesses multiple differentiability. Combined with the decay properties demonstrated in Theorem \ref{mainthm1}(iv), this enables us to apply a standard convergence result pertaining to the Gauss-Laguerre quadrature formula \citep{davis2007methods}. Consequently, we derive the desired result.
\end{pf}

\subsection{Implementation and Computational Details}\label{GLQuadratureimplimentation}
We are now ready to outline the method we propose for numerically computing $I^{\alpha}_a f(t_n)$, where $n = 1, 2, \cdots, P$. 
In this algorithm, the symbol $\hat{\phi}_l$ is utilized to represent the approximate value of $\hat{\phi}(x_l^{GLa},t_n)$ for the current time step, corresponding to the presently evaluated value of $n$. 

For given the initial point $a\in \mathbb{R}$, the order $\alpha\in (0,1)$, the grid points $t_n$, $n = 1, 2, \ldots, P$ and the number of quadrature nodes $\Lambda \in \mathbb{N}$,

\begin{enumerate}
  \item For $l = 1, 2, \ldots, \Lambda$:
  \begin{description}
    \item[a.] Compute the Gauss-Laguerre nodes $x_l^{GLa}$ and the associated weights $w_l^{GLa}$.
    \item[b.] Define the auxiliary terms $r_l \leftarrow \frac{-x_l^{GLa}}{1-\alpha}$ and $\tilde{r}_l \leftarrow \frac{x_l^{GLa}}{\alpha}$.
    \item[c.] Initialize $\phi_l \leftarrow 0$ and  $\tilde{\phi}_l \leftarrow 0$ to denote the initial condition for the differential equation (\ref{DE}).
  \end{description}
  \item For $n = 1, 2, \ldots, P$: 
   \begin{description}
    \item[a.] set $h \leftarrow t_n - t_{n-1}$.
    \item[b.] For $l= 1, 2, \ldots, \Lambda$:
     \begin{description}
    \item[(i)] Update the value of $\phi_l$ by solving the corresponding differential equation (\ref{DE}) using the backward Euler method
    \begin{equation}\label{algo1}
    \phi_l \leftarrow \frac{1}{1+h \mathrm e^{r_l}}[\phi_l + h c_{\alpha} \mathrm e^{(1-\alpha) r_l}f(t_n)].
    \end{equation}
    \item[(ii)] Similarly, update the value of $\tilde{\phi}_n$ by
       \begin{equation}\label{algo2}
    \tilde{\phi}_l \leftarrow \frac{1}{1+h \mathrm e^{\tilde{r}_l}}[\tilde{\phi}_l + h c_{\alpha} \mathrm e^{(1-\alpha)\tilde{r}_l} f(t_n)].
    \end{equation}
  \end{description}
  \item[c.]  Calculate the desired approximate value for $I^{\alpha}_a f(t_n)$ using the formula
     \begin{equation*}
   I^{\alpha}_a f(t_n) \approx \sum_{l=1}^{\Lambda} w_l^{GLa}\exp(x_l^{GLa})\bigg[\frac{1}{1-\alpha}\phi_l+\frac{1}{\alpha}\tilde{\phi}_l\bigg].
    \end{equation*}
  \end{description}
\end{enumerate}

Here we select the backward Euler method, considering the constant factor that multiplies the unknown function $\phi(\cdot,r)$ in (\ref{DE}), making it necessary to use an A-stable method. This choice is the simplest among such methods. Alternatively, one could opt for the trapezoidal method, another A-stable approach, necessitating adjustments to the formulas outlined in (\ref{algo1}) and (\ref{algo2}) as follows:
\begin{align}\label{algo3}
    \phi_l \leftarrow \frac{1}{1+\frac{h}{2} \mathrm e^{r_l}}&\bigg[\bigg(1-\frac{h}{2} \mathrm e^{r_l}\bigg)\phi_l\nonumber\\
    &+\frac{h}{2}c_{\alpha} \mathrm e^{(1-\alpha)r_l}[f(t_n)+f(t_{n-1})]\bigg],
\intertext{and}
\label{algo4}
    \tilde{\phi}_l \leftarrow \frac{1}{1+\frac{h}{2} \mathrm e^{\tilde{r}_l}}&\bigg[\bigg(1-\frac{h}{2} \mathrm e^{\tilde{r}_l}\bigg)\phi_l\nonumber\\
    &+\frac{h}{2}c_{\alpha} \mathrm e^{(1-\alpha)\tilde{r}_l}[f(t_n)+f(t_{n-1})]\bigg],
    \end{align}
respectively. Note that the effective evaluation of the formulas \eqref{algo1}--\eqref{algo4} in the given form may lead to overflows in some intermediate results. As shown in \citet{diethelm2021new}, this can be avoided by simple reformulations.
    
\section{Conclusion and Future Directions}
In this paper, our primary aim has been to craft highly efficient numerical techniques aimed at the approximate assessment of Riemann-Liouville fractional integrals with less computational complexity and memory footprint. For this, we have embarked on an exploration of innovative variations in diffusive representations tailored for fractional integrals. We have approximated the kernel function to the fractional integral by representing it as an exponential sum. This representation can be further optimized by leveraging Prony's method to curtail the number of terms involved. Subsequently, we have harnessed this refined approximation to compute an estimate for the fractional integral. This approach yields a notable reduction in both computational intricacy and memory usage, offering an enticing prospect for practical implementations. In addition to the exponential sum approximation, we have enriched our computational toolkit by developing the Gauss-Laguerre formula as an alternative method for approximating fractional integrals. 

This research opens several promising avenues for future investigations and applications such as extending the exploration of diffusive representations to other integral transforms and fractional operators to obtain valuable insights. Investigating more general forms of admissible transformations and their impact on the efficiency of numerical algorithms is an interesting direction. Furthermore, a comprehensive analysis of the error bounds and convergence properties of the numerical algorithms presented in this paper is essential. 

\bibliography{ifacconf}             % bib file to produce the bibliography
                                                     % with bibtex (preferred)
                                                   
%\begin{thebibliography}{xx}  % you can also add the bibliography by hand

%\bibitem[Able(1956)]{Abl:56}
%B.C. Able.
%\newblock Nucleic acid content of microscope.
%\newblock \emph{Nature}, 135:\penalty0 7--9, 1956.

%\bibitem[Able et~al.(1954)Able, Tagg, and Rush]{AbTaRu:54}
%B.C. Able, R.A. Tagg, and M.~Rush.
%\newblock Enzyme-catalyzed cellular transanimations.
%\newblock In A.F. Round, editor, \emph{Advances in Enzymology}, volume~2, pages
%  125--247. Academic Press, New York, 3rd edition, 1954.

%\bibitem[Keohane(1958)]{Keo:58}
%\newblock \emph{Power and Interdependence: World Politics in Transitions}.
%\newblock Little, Brown \& Co., Boston, 1958.

%\bibitem[Powers(1985)]{Pow:85}
%\newblock Is there a way out?
%\newblock \emph{Harpers}, pages 35--47, June 1985.

%\bibitem[Soukhanov(1992)]{Heritage:92}
%\newblock \emph{{The American Heritage. Dictionary of the American Language}}.
%\newblock Houghton Mifflin Company, 1992.

%\end{thebibliography}

%\appendix
%\section{A summary of Latin grammar}    % Each appendix must have a short title.
%\section{Some Latin vocabulary}              % Sections and subsections are supported  
                                                                         % in the appendices.
\end{document}